\documentclass[10pt]{amsart}       
\usepackage{amssymb}                  
\usepackage{amsthm}
\usepackage{graphicx}
\usepackage{amscd}
\usepackage{color}
\usepackage{enumerate}
\usepackage{xy}
\usepackage{amsmath}
\xyoption{all}

\theoremstyle{plain}
\newtheorem{theorem}{Theorem}
\newtheorem{lemma}{Lemma}
\newtheorem{corollary}{Corollary}

\theoremstyle{definition}
\newtheorem{definition}{Definition}

\newtheorem{example}{Example}
\theoremstyle{remark}

\newtheorem{remark}{Remark}
\def\ens{\ensuremath}                  
\newcommand\bydef{\ens{\stackrel{\text{def}}{=}}}
\newcommand{\N}{\mathbb{N}}
\newcommand{\Q}{\mathbb{Q}}
\newcommand{\R}{\mathbb{R}}

\newcommand{\Z}{\mathbb{Z}}

\newcommand{\HH}{\mathcal{H}}

\newcommand{\A}{\alpha}
\newcommand{\B}{\beta}

\newcommand\rb[2]{\raisebox{#1 mm}{#2}}                    
\newcommand\irrat{\ens{\R\rb{.26}{\text{$\setminus$}}\Q}}

\newcommand{\intpart}[1]{\big[ #1 \big]}                      
\newcommand{\intpartt}[1]{\big[ #1 \big]^{-}}                      

\newcommand{\fracpart}[1]{\langle#1 \rangle}                   
\newcommand{\card}{\mathbf{card}}                                
\newcommand{\modone}[1]{(\mathrm{mod}\ 1)}            
\newcommand{\cf}[1] {\big[ #1 \big]_\downarrow} 
\newcommand\uv[2]{\scalebox{#1}{#2}} 
\newcommand\uvm[2]{\scalebox{#1}{\ens{#2}}} 
\newcommand\hsp[1]{\mbox{}\hspace{#1mm}} 
\newcommand\vsp[1]{\par \vspace{#1mm}} 
\newcommand\vspm[1]{\par \vspace{-#1mm}} 
\newcommand\divides{\uv{1.2}{\,{\ens{\mid}}}\,} 
\newcommand\miduv[1]{\uv{1.#1}{\,{\ens{\mid}}}\,} 
\newcommand\uoi{\mathbb{I}_{_{0,1}}\!}  
\def\bsq{\hfill\raisebox{-3mm}{$\blacksquare$}\;}

\newcommand\ppp[1]{}  
\title{Continued Fractions and Heavy Sequences}

\author[M.\ Boshernitzan]{Michael Boshernitzan}

\address{Department of Mathematics, Rice University, Houston, TX~77005, USA}

\email{michael@rice.edu}

\author[D.\ Ralston]{David Ralston}

\address{Department of Mathematics, Rice University, Houston, TX~77005, USA}

\email{dsr@rice.edu}

\thanks{D.\ R.\ was supported in part by NSF VIGRE grant
DMS--0240058.}

\begin{document}
\begin{abstract} 

We initiate the study of the sets   $\HH(c)$, $0\!<\!c\!<\!1$,   of  real  $x$   for which 
the sequence   $(kx)_{k\geq1}$  (viewed mod~1)  consistently hits 
the interval  $[0,c)$  at least as often as expected (i.\,e., with frequency  $\geq c$).
More formally,
$$
\HH(c)\bydef\big\{\A\in \R  \miduv2  \card\big(\{1\leq k\leq n\mid 
         \fracpart{k\A}<c\}\big)\geq cn, \ \text{ for all } n\geq1\big\},
$$
where  $\fracpart{x}=x-\intpart{x}$  stands for the fractional part of  $x\in\R$.

We prove that,  for rational  $c$,  the sets $\HH(c)$  are of positive 
Hausdorff dimension and,  in particular,  are uncountable. 
For  integers $m\geq1$,  we obtain a surprising
characterization of the numbers  $\A\in\HH_m\bydef\HH(\tfrac1m)$  in terms
of their continued fraction expansions:  The odd entries  (partial quotients)   
of these expansions  are divisible by~$m$.
The characterization  implies that $x\in\HH_m$    if and only if  
$\frac 1{mx} \in\HH_m$,  for $x>0$.
We are unaware of a direct proof of this equivalence, without making a use 
of  the mentioned characterization of  the sets $\HH_m$.

We also introduce the dual sets   $\widehat\HH_m$  of reals  $y$  for which 
the sequence  of integers  $\big([ky]\big)_{k\geq1}$  consistently hits the 
set  $m\Z$  with the at least expected frequency $\frac1m$  and establish
the connection with the sets  $\HH_m$:
$$ \text{
If \ $xy=m$  \ for \  $x,y>0$, \  then  $x\in\HH_m\iff y\in\widehat\HH_m$.}
$$
The motivation  for the present study comes  from Y.~Peres's ergodic lemma in \cite{peres}.
 \end{abstract}
 

\maketitle
                                                                                                                                                                             \ppp{Section 1\\ Defs 1}
\section{Notation and Results}

We write \ $\R\supset\Q\supset\Z\supset\N$ \  for the sets of real numbers,  rational numbers,
integers and positive integers respectively. \vsp1

In the paper we initiate the study of the sets   $\HH(c)$, $0\!<\!c\!<\!1$,   
of  $x\in\R$   for which the sequence   $(kx)_{k\geq1}$  (viewed mod~1)  consistently 
hits the interval  $[0,c)$  at least as often as expected.  More formally,
                                                                                                                               \ppp{eq:hc - 2} 
\begin{equation}\label{eq:hc}
\HH(c)=\big\{\A\in \R \ \miduv2 \  \card\big(\{1\leq k\leq n\mid \fracpart{k\A}<c\}\big)
             \geq cn, \ \text{ for all } \ n\in\N\big\}
\end{equation}
where  $\fracpart{x}=x-\intpart{x}$  stands for the fractional part of  $x\in\R$.
Define
\begin{equation}\label{eq:hhm}
\HH_m=\HH(\tfrac1m), \; \text{ for }  m\in\N.
\end{equation}

The following notations will be used for  CF  (continued fraction) expansions                                               \ppp{CFs - 3}
of finite length \  $n+1$:
$$
  \cf{a_0,a_1,a_2,\ldots,a_{n}}=a_0+\frac{1}{a_1+\frac{1}{a_2+\frac{1}
        {\ddots \,+ \frac{1}{a_n}}}}; \quad n\geq0,
$$
\vspm3
\noindent  or of infinite length 
$$
 \cf{a_0,a_1,a_2,\ldots}=\lim_{n \rightarrow \infty} \cf{a_0,a_1,a_2,\ldots,a_n},
$$
where  $a_0 \in \Z$ and $a_k\in\N$  \ for  $k\geq1$.

For some basic facts and standard notation from the theory of  CFs  we refer to                                         \
\cite{lang}  or  \cite{khinchin}.  (First few pages in  either book  should suffice
for our purposes.)

Every irrational number has a unique infinite CF expansion, and every 
rational number has exactly two finite  CF expansions
$$
\cf{a_0,a_1,a_2,\ldots,a_{n-1}+1}=\cf{a_0,a_1,a_2,\ldots,a_{n-1},1},\quad
  \text{with  }a_{n-1} \geq 1, 
$$
(with the lengths being two consecutive integers, $n$ and $n+1$).                            \ppp{OCF\\ECF- 5}

\begin{definition}  \label{def:ocf}
By the  {\em odd}\, CF  (odd continued fraction)  expansion  (of  $\A\in\R$)   
we mean   the  CF   expansion  of   length \  $L\in\{\infty,1,3,5,\ldots\}$.   
Similarly,  in the   {\em even} CF  expansions
one  assumes   $L\in\{\infty,2,4,6,\ldots\}$. 
\end{definition}

This way every number $\A\in\R$  has  unique  
both odd and even CF  expansions;
the two coincide if and only if  $\A$  
is irrational.  The sequence of  (CF)  convergents of\, $\A$                                                    \ppp{convergents \\ 6}
$$
\delta_k(\A)=\cf{a_0,a_1, \ldots,a_k}, \quad 0\leq k<L,
$$
can be  alternatively  defined  as the sequence of rational numbers \  $\delta_k=\frac{p_k}{q_k}$ \
with  numerators and denominators  $p_k=p_k(\A),\, q_k=q_k(\A)$  
determined  by the recurrence relations                								       \ppp{recurr. rels\\   eq:RR \\ 7}		
\begin{equation} \label{eq:RR}
\begin{cases}
p_k=a_kp_{k-1}+p_{k-2}, \\
q_k=a_kq_{k-1}+q_{k-2},  
\end{cases}
\; \text{for \ }  2\leq k<L,
\end{equation}
and the initial conditions  \;  
$   p_0=a_0; \; q_0=1;  \quad p_1=a_0 a_1+1; \;  q_1=a_1.  $

The following theorem provides a criterion for the relation  $\A\in\HH_m$  
to hold (see \eqref{eq:hhm}).                                                                                               \ppp{Criterion \\ thm:hhm \\  8}

\begin{theorem} \label{thm:hhm}
Let   $\A\in\R$  and  assume that  $\A=\cf{a_0,a_1,a_2,\ldots}$  is its odd CF                                                    
expansion  (i.\,e., of the length  $L\in\{\infty,1,3,5,\ldots\}$).  Let  $m\in\N$  be given.
Then the following three conditions are equivalent: \vsp1

\hsp{5} \parbox{150mm}{
\begin{itemize}
\item[\bf (C1)]   $\A\in\HH_m$. \vsp1
\item[\bf (C2)]    $m\divides a_{k}$,  for all odd  $k$,  $1\leq k < L$. \vsp1
\item[\bf (C3)]    $m\divides q_{k}$,  for all odd  $k$,  $1\leq k < L$,   where  $q_{k}=q_{k}(\A)$  
          are the denominators of the convergents for\,  $\A$,  see \eqref{eq:RR}.
\end{itemize}
}
\end{theorem}


\begin{remark}  For  $m=1$  the above theorem holds  trivially  because   $\HH_1=\R$.                            \ppp{Rem 1 \\ 9}
It also holds trivially for  $\A\in\Z$  (in this case  $L=1$).
\end{remark}
{\bf Examples.}  \begin{enumerate}
\item  $\A=\tfrac43$.  The odd CF expansion is  $\cf{1,2,1}$, $L=3$, $a_1=2$.  
             Thus  $\tfrac43\in\HH_m$ if and only if  $m=1\text{ or } 2$.
\item  $\A=\frac{\sqrt{5}}2$.  The odd CF expansion is  $\cf{1,8,2,8,2,8,\ldots}$, $L=\infty, \ a_1=a_3=a_5=\ldots=8$.  \\
            Thus  $\frac{\sqrt{5}}2\in\HH_m$\, if and only if\, $m=1, 2, 4 \text{ or } 8$.
\end{enumerate}

\begin{corollary}\label{cor:1}
$\HH_m\cap\HH_n=\HH_{\text{\rm{LCM}}(m,n)}$, \ for all  \  $m,n\in\N$.                                               \ppp{cor:1 \\ 10}
\end{corollary}

\begin{corollary}\label{cor:2}
For  real $\A>0$ and   $m\in\N$,  we have  \  $\A\in\HH_m$ if and only 
if \  $\frac 1{m\A}\in\HH_m$.                    \ppp{cor:2 \\ 11}
\end{corollary}

Both corollaries follow directly from the equivalence of  {\bf (C1)}  and    {\bf (C2)}   in
Theorem \ref{thm:hhm};  the proof of Corollary~\ref{cor:2}  
also uses the identity                                                                                                                                           \ppp{eq:id \\ 12}
\begin{equation}\label{eq:id}
m\cf{x_0, mx_1, x_2, mx_3, x_4, mx_5, \ldots}=\cf{mx_0, x_1, mx_2, x_3, mx_4, x_5, \ldots}.
\end{equation}

In the next three theorems  we classify  the numbers 
in the sets   $\widehat\HH_m$,  $m\in\N$:
\begin{equation}\label{eq:hbarc}
\widehat\HH_m= \big\{\A\in \R \ \miduv2 \  \card\big(\{1\leq k\leq n\mid \intpart{k\A}\in 
         m\Z\}\big)\geq\tfrac nm,  \ \text{ for all  } n\in\N\big\}.
\end{equation}
\vsp1

\begin{theorem}\label{thm:hbarh}
For   $\A\in\R$ and\, $m\in\N$,    we have  \  $\A\in\HH_m \iff m\A\in\widehat\HH_m$.
\end{theorem}

The proof of Theorem \ref{thm:hbarh} is derived from the comparison \eqref{eq:hc}                \ppp{Thm 2 \\thm:hbarh\\ 13}
and  \eqref{eq:hbarc}  and taking in account that,  for  $x\in\R$,                                                                 
$
\fracpart{x}\in [0,1/m) \iff \intpart{mx}\in m\Z.
$

Note that we establish another,  deeper connection 
(than the one indicated in Theorem \ref{thm:hbarh})  between the sets  
$\HH_m$  and $\widehat\HH_m$  in  Theorem \ref{thm:hbarh2}  below.                                          \ppp{next\\thm \\ barhm\\14}

The following result provides explicit description of  the sets  $\widehat\HH_m$.          

\begin{theorem} \label{thm:barhm}                                                         
Let   $m\in\N$,  $\A\in\R$  and  assume that  $\A=\cf{a_0,a_1,a_2,\ldots}$  is  its even CF                  \ppp{thm\\ barhm\\15}                                                   
expansion (of the length  $L\in\{\infty,2,4,6,\ldots\}$).  Let  .
Then the following three conditions are equivalent: \vsp1

\hsp{20} \parbox{115mm}{
\begin{itemize}
\item[\bf (C1)]   $\A\in\widehat\HH_m$.   \vsp1
\item[\bf (C2)]    $m\divides a_{k}$,  for all even  $k, \ 0\leq k < L$. \vsp1
\item[\bf (C3)]    $m\divides p_{k}$,  for all even  $k, \ 0\leq k < L$
            where  $p_{k}=p_{k}(\A)$  are numerators of the convergents for  $\A$,  
            see \eqref{eq:RR}.
\end{itemize}
}
\end{theorem}

The proof of Theorem \ref{thm:barhm}  easily follows from Theorems \ref{thm:hhm} 
and \ref{thm:hbarh}  using the identity   \eqref{eq:id}.

Alternatively,  Theorem \ref{thm:barhm}  can be derived from the following                              \ppp{thm 4\\ hbarh  \\ 16}

\begin{theorem}\label{thm:hbarh2}
For   $\A>0$  and  $m\in\N$,    we have  \  $\A\in\HH_m \iff \frac1\A\in\widehat\HH_m$.
\end{theorem}

Theorem \ref{thm:hbarh2}  follows from Corollary \ref{cor:2}  and identity \eqref{eq:id}. 
\vsp2

The proof of Theorem \ref{thm:hhm}   will be provided in the next section.  We also prove 
(Theorem \ref{thm:suf} and \ref{thm:hhm2})  that
$$ \HH(\tfrac nm)\supset  \HH(\tfrac 1m)=\HH_m, \; \text{ for arbitrary } \;n,m\in\N, \; n<m,
$$
and conclude that, for rational\,  $c,\ 0<c<1,$  the sets  $\HH(c)$  have a positive Hausdorff dimension (Corollary~\ref{cor:3}).

Finally, in the last section we discuss briefly the motivation behind our study. \vsp3

\section{Proof of Theorem  \ref{thm:hhm}.}

The proof is subdivided into several lemmas,  some of which are of  independent interest.                 \ppp{eq\\sna\\17}
Let\,    $\uoi$\,  stand for the open  unit interval\,  $(0,1)$.
For\,  $n\in\N$, $\A>0$\,  and\,  $c\in\uoi$,
consider the following  finite subsets of\,  $\N$:                                                                                                                    
\begin{equation}\label{eq:sna}
   S(n,\A)\bydef\{k\in \N\miduv1 k\A<n\}
\end{equation}
and                                                                                                                                                                  \ppp{eq\\snac\\18}
\begin{equation}\label{eq:snac}
   S(n,\A, c)\bydef\{k\in S(n,\A)\miduv1 \fracpart{k\A}<c\}=\{k\in \N\miduv1 k\A<n
         \ \& \ \fracpart{k\A}<c\}.
\end{equation}

It is easy to see that                                                                                                                                    \ppp{eq\\cardsna\\19}
\begin{equation}\label{eq:cardsna}
   \card\big(S(n,\A)\big)=\uv{1.1}{\ens{\intpartt{\tfrac n\A}}}
\end{equation}
and                                                                                                                                                               \ppp{eq\\cardsnac\\20}
\begin{equation}\label{eq:cardsnac}
     \card\big(S(n,\A,c)\big)=\uv{1.1}{\ens{\intpartt{\tfrac c\A}}}+\sum_{k=1}^{n-1}
     \uv{1.1}{\ens{\big(\intpartt{\tfrac{k+c}\A}-\intpartt{\tfrac k\A}\big)}}
\end{equation}
where  $\intpartt{x}$  stands for the largest integer  smaller than \ $x\in\R$:                                     \ppp{eq\\intpartt\\21}
\begin{equation}\label{eq:intpartt}
\intpartt{x}\ \bydef  \
\begin{cases}
\intpart{x} \qquad \text{ if } x\notin\Z.\\  
x-1  \quad   \text{ if } x\in\Z.
\end{cases}
\end{equation}

We observe the following                                                                                                                     \ppp{Lemma\\lem:1\\ 22}

\begin{lemma}\label{lem:1}
Given  \   $\A>0$  \  and \  $c\in\uoi=(0,1)$,  the following two conditions are equivalent: 

\hsp{12} \parbox{110mm}{
\begin{enumerate}
\item[\rm(1)]  $\A\in\HH(c)$;
\item[\rm(2)]  $\card\big(S(n,\A,c)\big) \geq  c \ \card\big(S(n,\A)\big)$,  \  for all  \ $n\in\N$.
\end{enumerate}
}
\end{lemma}
\noindent {\bf Proof}. The claim of Lemma \ref{lem:1}  follows directly from the definitions             \ppp{Proof\\ lem:1\\ 23}
of the sets  \  $\HH(c)$, $S(n,\A)$  and  $S(n,\A,c)$  (see   \eqref{eq:hc},  
\eqref{eq:sna},  \eqref{eq:snac}).                          \bsq

\begin{lemma}\label{lem:2}
Let  $\A,\B>0$  and  $c\in\uoi$.  Assume that the following two conditions are met:                             \ppp{eq:rat \\ 24}
\begin{equation}\label{eq:rat}
\text{\rm \small(1)}\ \tfrac1\A-\tfrac1\B\in \Z \qquad \text{  and } \qquad  
          \text{\rm \small(2)}\ \tfrac c\A-\tfrac c\B\in \Z.
\end{equation}  
\noindent Then:

\hsp{12} \parbox{110mm}{
\begin{itemize}
\item[\rm (A)] \hsp1 $\card\big(S(n,\A)\big)-\card\big(S(n,\B)\big)=\frac{n(\B-\A)}{\A\B}${\rm ;}
\item[\rm (B)] \hsp1 $\card\big(S(n,\A,c)\big)-\card\big(S(n,\B,c)\big)=c\,\frac{n(\B-\A)}{\A\B}${\rm ;}
\item[\rm (C)] \hsp2 $\A\in\HH(c)\iff\B\in\HH(c)$.
\end{itemize}
}
\end{lemma}

\noindent{\bf Proof}.  The claims   (A) and  (B)  of the  lemma follow                                                 \ppp{proof \\ lem:2 \\ 25}
from formulae  \eqref{eq:cardsna}, \eqref{eq:cardsnac}  and the obvious implications
$$x-y\in \Z\ \iff\  \fracpart{x}= \fracpart{y} \ \implies \ \intpartt{x}\!-\intpartt{y}\!=x-y.$$

\noindent Finally,  (C)  follows from  (A), (B)  and Lemma \ref{lem:1}.    \bsq
\vsp3

The next lemma is just a specification  of Lemma \ref{lem:2}.                                                           \ppp{speci-\\fication \\ 26}
\vsp1

\begin{lemma}\label{lem:3}
Let  \  $\A,\B>0$, $c\in\uoi$,  $m\in\N$  \  be given
and assume that the following two conditions are met:
$$
\text{\rm (a)}\ \tfrac1\A-\tfrac1\B\in m\,\Z, \hsp{12} \text{\rm (b)}\ mc\in\N. \hsp{25}
$$

Then \ $\A\in\HH(c)$\, if and only if\, $\B\in\HH(c)$.
\end{lemma} 
\noindent{\bf Proof}.  
We need only  validate  condition  (2) of    Lemma \ref{lem:2}.                                                                         \ppp{proof\\lem:3\\27}
It follows from   (a)   that   \  $\tfrac1\A-\tfrac1\B=km$, for some  $k\in\Z$.
But then   $\tfrac c\A-\tfrac c\B=ckm=k(mc)\in\Z$,  in view of  (b).  
The proof is complete.     \bsq
\vsp3

For  $\A\in\R$,  denote by  $L(\A)$  the length of the odd CF expansion of  \  
$\A=\cf{a_0(\A),a_1(\A),\ldots}$.
Observe that $L(\A)=1$ \ if and only if  \ $\A\in\Z$;  otherwise   $L(\A)\geq3$.
Define the map  \  $\phi:  \R\to\R$  \  by the rule
$$
\phi(\A)=
\begin{cases}
\frac1{\fracpart{\A}}=\cf{a_1,a_2,\ldots}, \hsp{2} \text{if }\ \A\notin\Z, \\ 
\hsp2 0,  \hsp{27} \text{if }\ \A\in\Z.
\end{cases}
$$

One easily verifies that for  $\A=\cf{a_0, a_1, a_2,\dots}\notin\Z$  \  one has  \
$
\phi^2(\A)=\cf{a_2,\dots}
$,  i.\,e.   $\phi^2$ \  removes the first two entries (partial quotients) 
in the odd CF  expansion of any non-integer.
\vsp2

Next   we introduce the sets                                                                                                                   \ppp{eq:rm\\28}
\begin{equation}\label{eq:rm}
\R(m)\bydef\big\{\A\in\R\!\setminus\!\Z\,\miduv2\,a_1(\A)\in m\,\N\big\}\uvm{1.5}\cup\,\Z, 
\quad  m\in\N.
\end{equation}

\begin{lemma}\label{lem:4} 
Let   $m\in\N$,  $c\in\uoi$   and assume
that \  $mc\in\N$.  Then,  for every \  $\A\in\R(m)$, 
$$\A\in\HH(c)\iff\phi^2(\A)\in\HH(c).$$
\end{lemma} 

\noindent{\bf Proof}.                                                                                                                                 \ppp{proof\\lem:4\\29}
Denote  $\B=\phi^2(\A)$  and  $u=\tfrac1\A-\tfrac1{\fracpart{\B}}=-a_1$.  
Since  $\A\in\R(m)$,  we conclude that  $m\!\divides\! u$  and use  Lemma~\ref{lem:3}
to complete the proof of Lemma \ref{lem:4}:  \;
$
\A\in\HH(c)\iff \fracpart{\B}\in\HH(c)\iff \B\in\HH(c). $
\bsq         
\vsp2

It turns out that Lemma \ref{lem:4} can be used  to explicitly exhibit
uncountable subsets of   $\HH(c)$,  for $c\in \uoi\cap\Q$.
($\Q$  stands for the set of rational numbers.)
Those are the sets
\begin{equation}\label{eq:rm2}
\R_m\bydef\{\A\in\R \miduv1 \phi^{2k}(\A)\in \R(m),  \;  \text{ for all } \ k\geq0\}, 
\quad   m\in\N,
\end{equation}
(see  Theorem \ref{thm:suf}   below).  \vsp1

The following lemma provides an alternative, 
more explicit description of the classes\,  $\R_m$.                                                                                   \ppp{lem:5\\30}

\begin{lemma}\label{lem:5} 
Let  $m\in\N$  and assume that  \  $\A\in\R$  is given in terms of its odd  CF  expansion 
$\A=\cf{a_0,a_1,\ldots}$ of length \text{$L=L(\A)\in\{\infty,1,3,\ldots\}$.}  Then 
$\A\in\R_m$  if and only if  $a_k\in m\Z, \;  \text{\rm for all odd }  k<L.
$
\end{lemma} 

\noindent{\bf Proof}.                                                                                                                                       \ppp{proof\\lem:5 - 31}
The proof follows directly from the nature of the map  $\phi^2$  and the trivial fact 
that  $\Z\subset\R(m)$.                              \bsq  
\vsp2

\noindent{\bf Examples.}  
\begin{enumerate}
\item  $\A=\tfrac43$.  The odd CF expansion is  $\cf{1,2,1}$, $L=3$, $a_1=2$.  
          Thus  $\tfrac43\in\R_m \iff  m=1\text{ or } 2$.
\item  $\A=\frac{\sqrt{5}}2$.  The odd CF expansion is  
           $\cf{1,8,2,8,2,8,\ldots}$, $L=\infty, \ a_1=a_3=a_5=\ldots=8$.  \\
           Thus  $\frac{\sqrt{5}}2\in\R_m \iff  m=1, 2, 4 \text{ or } 8$.
\end{enumerate}
\vsp1

\begin{theorem}\label{thm:suf}
Let $m\in\N$,  $c\in\uoi\!=(0,1)$   and assume
that \  $cm\in\N$.  Then \
$\R_m\subset \HH(c).$
\end{theorem}

\noindent{\bf Proof of Theorem \ref{thm:suf}}.   
Recall that  $\Q\subset\R$  stands for the set of rational numbers.  
We prove that
$$\A\in\R_m\implies\A\in\HH(c).$$                                                                                                                   \ppp{proof\\thm 5\\suf - 32}

\noindent{\bf Case 1}.  $L<\infty$, i.\,e.  $\A\in\Q$. The proof goes 
by induction in\,  $L=L(\A)\in\{1,3,5,\ldots\}$.    

If  $L=1$, then  $\A\in \Z$  and one has
both  $\A\in\R_m$  and  $\A\in \HH(c)$.  
For the inductional step, we use  Lemma \ref{lem:4} and the obvious
fact that   $\phi^2(\R_m)\subset\R_m$   (see \eqref{eq:rm2} and Lemma \ref{lem:5}).

\noindent{\bf Case 2}.  $L=\infty$, i.\,e.  $\A\in\R\!\setminus\!\Q$\,  ($\A$\, is irrational).  The proof uses 
an approximation argument.  For a subset  $S\subset\R$,  denote 
by\,  $\overline S$\,    the closure of  $S$  in  $\R$. 
Next  we validate  the following inclusion:   
\begin{equation}\label{eq:valid}
\overline{\HH(c)}\cap(\irrat)\subset\HH(c)  \qquad (\text{for }\,  c\in\Q\cap(0,1)).
\end{equation}
Indeed, it follows from \eqref{eq:hc}  that
$
\HH(c)=\bigcap_{n\in\N}\, \HH(c,n), 
$
where   
$$
\HH(c,n)=\big\{\A\in\R\mid\card(\{1\leq m\leq n\mid\fracpart{m\A}<c\})\geq cn\big\}
$$
is a finite union of intervals of the form\,  $[u, v)$\,   with rational endpoints\,  $u, v$:
$$
u, v\in\!\!\bigcup_{0\leq r< s\leq n}\left\{\tfrac rs, \tfrac {r+c}s\right\}\subset\Q
$$
because\,    $c\in\Q$.
In particular,   for all  $n\geq1$, 
$$
\overline{\HH(c,n)}\cap(\irrat)\subset\HH(c,n).
$$

The proof of    \eqref{eq:valid}  is completed as follows:
$$
\overline{\HH(c)}\cap(\irrat)=\overline{\bigcap_{n\in\N}\, \HH(c,n)}\cap(\irrat)\subset
\Big(\bigcap_{n\in\N}\overline{\HH(c,n)}\,\Big)\cap(\irrat)=
$$
$$
\hsp{15}=\bigcap_{n\in\N}\left(\overline{\HH(c,n)}\cap(\irrat)\right)
\subset\bigcap_{n\in\N}{\HH(c,n)}=\HH(c).
$$

Next one considers the sequence\,  $\delta_{2k}=\delta_{2k}(\A),\, k\geq0,$\,   of even  CF  
convergents  of\,  $\A$\,  (with   $L(\delta_{2k})=2k+1$, an odd number).  
By what has been proven in Case 1,  \  $\delta_{2k}\in\R_m\cap\,\Q\subset\HH(c)$, 
for all\,  $k\geq0$. 
Since\,   $\lim_{k\to\infty} \delta_{2k}=\A\in\irrat$, the proof is complete in view of  \eqref{eq:valid}.								      	         	   \bsq  

\begin{corollary}\label{cor:3}
Let  \ $C\subset\Q\cap\uoi$  be a finite subset of rational numbers.  Then the intersection \                                                          \ppp{cor:3\\?? - 33}
  $\bigcap_{c\in C} \HH(c)$  is an uncountable   set of positive Hausdorff dimension.
\end{corollary}

\noindent{\bf Proof}. Let  $m\in \N$  be a common denominator 
for all the numbers  $c\in C$.  Then $\R_m\subset\bigcap_{c\in C} \HH(c)$, 
in view of Theorem \ref{thm:suf}.                                                                                                                              \ppp{reference\\?? - 34}

The set  $\R_m$  is clearly uncountable,  and it is easily
seen  to have a positive Hausdorff dimension.   
(One way to see it is to observe that it contains the set
$\R'_m$   of all numbers of the form
$\cf{0,m,a_1,m,a_2,m,\ldots}$  where   $a_i \in \{1,2\}$.
This set is the disjoint union of its images under the
two contractions: 
$$f_i(x)=\frac{1}{m+\frac{1}{i+x}}, \hsp3 i=1,2,$$
and therefore must be of positive Hausdorff dimension
(see Chapter 9 in \cite{falconer}).)
    \bsq  

\vsp1

\begin{corollary}\label{cor:rmhm}
$\R_m\subset\HH_m$,  for all \ $m\in\N$.                                                                                                                \ppp{cor\\rmhm - 35}
\end{corollary}

\noindent{\bf Proof of Corollary  \ref{cor:rmhm}}.                                                                                                                                                   
Take   $c=\tfrac1m$  in   Theorem \ref{thm:suf}.  
(Recall that  $\HH_m=\HH(\frac1m)$.)         	         	         	         	       \bsq  
\vsp2

As it is pointed out in Section 1,  the inclusion in Corollary \ref{cor:rmhm}  
can be reversed.                                                                                                                                                        \ppp{thm\\hhm2 - 36}

\begin{theorem}\label{thm:hhm2}
$\R_m=\HH_m$,  for all \ $m\in\N$.
\end{theorem}

We first need to prove the following

\begin{lemma}\label{lem:6}
$\HH_m\subset\R(m)$,  for all \ $m\in\N$.                                                                                                                 \ppp{lem:6 - 37}
\end{lemma}

\noindent{\bf Proof of Lemma \ref{lem:6}}. Since  $\Z\subset\R(m)$  (see  \eqref{eq:rm}),  
it suffices to prove that if  \ $\A\in\!\big(\R\!\setminus\!\Z\big)\cap\,\HH_m$  \  
then \  $\A\in\R(m)$.   We assume  without loss of generality that  $\intpart{\A}=0$  
(otherwise replacing  $\A$  by  $\fracpart{\A}$).
Let   $\cf{0,a_1,a_2,\ldots}$  be  the odd  CF  expansion of  $\A$,  we have to show 
that  $m\divides a_1$.

If not,  then  $a_1\equiv r \pmod m$, for some integer \ $r$, $1\leq r\leq m-1$.  
Define  $\B\in\R$  by its  CF  expansion  $\cf{0,r,a_2,a_3,\ldots}$,  with all the entries  $a_k$, 
for  $k\neq1$,  being  the same as in the  CF  expansion of  \ $\A$.  

Then, with  $c=\frac1m$,  we easily validate the conditions of Lemma \ref{lem:3}.
Indeed,  $\A\in\HH_m=\HH(c)$,  $mc=1\in\N$  \ and \  $\frac1\A-\frac1\B=a_1-r\in m\Z$ . 
By Lemma \ref{lem:3},  $\B\in\HH(c)=\HH_m$,  which is impossible because
$\B=\fracpart \B>\frac1{r+1}\geq\frac1m=c,$  so that  the relation   $\B\in\HH(c)$                  \ppp{lem:6\\proof\\end - 37}
contradicts  \eqref{eq:hc}   for  $n=1$.   \bsq
\vsp4

\noindent{\bf Proof of Theorem \ref{thm:hhm2}}.  
In view of Corollary \ref{cor:rmhm},  we only have to establish the inclusion  
$\HH_m\subset\R_m$.  Let  $\A\in\HH_m$  be given. We claim that then
\begin{equation}\label{eq:f2k}
\phi^{2k}(\A)\in\R(m)\cap\HH_m, \; \text{ for all } \  k\geq0.
\end{equation}
The proof goes by induction in $k$.  For  $k=0$,  \eqref{eq:f2k}
holds in view of Lemma \ref{lem:6}.  For the inductional step,
we use  Lemmas  \ref{lem:4}  and  \ref{lem:6}.
This completes the proof of  \eqref{eq:f2k}.
In view of  \eqref{eq:rm2},  $\A\in\R_m$, completing the proof.  
\bsq  
\vsp3

Now we are ready to complete the proof of Theorem  \ref{thm:hhm}   in the introduction.  
The main part of work  has been already done:  Theorem  \ref{thm:hhm2}  is a rephrasing 
of the equivalence  {\bf (C1)$\iff$(C2)}.

It remains to prove the equivalence  of the following two conditions:

\hsp6 {\bf (C2)} $m\divides a_k$,  \ for all odd \  $k$, $1\leq k< L$.

\hsp6 {\bf (C3)} $m\divides q_k$,  \ for all odd \  $k$, $1\leq k< L$.
\vsp1

The proof goes by induction in  $k$.  For  $k=1$ the equivalence is immediate 
because  $q_1=a_1$. 

Now assume that both  {\bf (C2)} and {\bf (C3)}   hold for some odd\, $k=n<L-2$.  
It suffices to show that
$$
m\divides a_{n+2} \iff m\divides q_{n+2}.
$$
\indent From the identity\,    $q_{n+2}=a_{n+2}q_{n+1}+q_n$\,  we derive the congruence\,  
$q_{n+2}\equiv a_{n+2}q_{n+1}\! \pmod m$,  so that the implication  $\implies$  is immediate.  
The opposite implication is also valid because\,  $q_{n+2}, q_{n+1}$\,  are relatively prime.
\vsp5

\section{Motivation: Heavy sequences.}

The following result by Y. Peres is closely related to the Maximal Ergodic Theorem:
\begin{lemma}[Peres]\label{lem:peres}
Let $T:X \rightarrow X$ be a continuous transformation of a compact space, and let $\mu$ be a probability measure preserved by $T$.  For every continuous $g:X \rightarrow \R$ there exists some $x \in X$ such that $$\forall N \in \N \hspace{.2 in}\frac{1}{N}\sum_{n=0}^{N-1}g(T^n x) \geq \int_X g d\mu.$$
\end{lemma}

This lemma can be strengthened to upper semi-continuous functions (such as characteristic functions of closed sets), and has several interesting generalizations - see \cite{ralston}.  The application of most interest to us right now is to extend the ideas of equation   \eqref{eq:hc}   to the context of dynamical systems.  Let $\{X,\mu\}$ be a compact space with probability Borel measure $\mu$, and let $T$ be a continuous map from $X$ to itself which preserves the measure $\mu$: $\mu(T^{-1}A)=\mu(A)$ for any Borel set $A$.  Define the \emph{heavy set of $f$}:

\begin{align}
\mathcal{H}_T^f = \big\{x \in X \textrm{: } S_n(x)-n \int_X f d\mu \geq 0,  \hsp2 \forall n \in \mathbb{N}\big\}.
\end{align}

Then Lemma \ref{lem:peres} tells us that in this situation, for any $f \in L^1(X,\mu)$, $\mathcal{H}_T^f \neq \emptyset$.  We also say that there is some point $x$ whose orbit is \emph{heavy for $f$}.  
  If $f$ is the characteristic function of a set $A$, we will generally simply refer to ``the heavy set of $A$," or call a sequence ``heavy for $A$."  Restricting ourselves only to the reals modulo one, $\R/\Z=S^1$, we derive the following results:

\begin{example}{Fix $\alpha \in S^1$.  Then for any closed subset $A\subset S^1$, there exists some point $x \in S^1$ such that the sequence $\{x+n\alpha\}_{n=0}^{\infty}$ is heavy for $A$.}
\end{example}

The previous example can be viewed as the following:  for any choice of a closed set $A\subset S^1$ and leading coefficient $\alpha$, there exists some choice of $\beta$ such that the polynomial $\alpha n + \beta$, considered modulo one, is heavy for $A$.
 This example may be generalized as follows:

\begin{example}{Fix $\alpha \in \mathbb{R}$, a closed set $A \subset S^1$, and a choice of $k \in \mathbb{N}$.  Then there exists a choice of coefficients $a_0, a_1, \ldots,a_{k-1}$ such that the sequence $$\{\alpha n^k + a_{k-1}n^{k-1} +\ldots +a_1 n + a_0\}_{n=0}^{\infty}$$ is heavy for $A$ (when taken modulo one).}
For details on how to derive this sequence as the orbit of a measure preserving system, we refer the reader to  p.Ê35--37  in  \cite{furstenberg},  or,  for a more detailed derivation of heaviness properties, to \cite{ralston}.
\end{example}

Finally, the reader may be tempted to try to generalize the results of Theorem \ref{thm:hhm} and
Example 1 to claim that the set
$$
\HH[A]\bydef \{x\in S^1 \miduv1 \text{ the sequence } (kx)_{k\geq1}  \ \text{ is heavy for }  A\}
$$
is always non-empty. This cannot be done.

\begin{example}\label{ex:A}
There exists a closed set $A \subset S^1$, a finite union of closed intervals, whose measure is larger than $1/3$, such that \; $\left(x \in A\right) \Rightarrow \left( 2x \notin A \textrm{, }3x \notin A\right).$
In particular,  for such  $A$  one has   $\HH[A]=\emptyset$.
\end{example}

The measure of the set  $A$  in the above example can be made arbitrary
close to  $\frac12$.   For details, see \cite{ralston}  where techniques from ergodic theory
are used (in a noncostructive way)  to establish the existence of such a set  $A$.

We have  examples of closed subintervals   $J\subset (0,1)\subset S^1$   for which the set
$\HH[J]$   is countable  or  finite   ($J=[\frac13, \frac23]$   and   $J=[\frac27, \frac57]$,
respectively).  We don't know whether it can be made empty.

Note that the subject of our paper is somehow related to that 
in \cite{dupain-sos}  where some sufficient conditions 
for the one-sided boundedness of the sequence
$$\card(\{1\leq k\leq n\miduv1\fracpart{k\A}<c\})-cn$$
have been established  (cf.~equation \eqref{eq:hc}).  
All our results are new and imply some of the results in  \cite{dupain-sos}.

\end{document}